\documentclass{article}

\usepackage{graphicx}      
\usepackage{amsmath}
\usepackage{amssymb}

\title{Modeling and Stabilization of a Rotating Mechanical System with Elastic Plates}

\author{Alexander Zuyev$^{1,3}$ and Julia Novikova$^{2,3}$}
\date{\small $^{1}$Max Planck Institute for Dynamics of Complex Technical Systems,\\
   Magdeburg, Germany\\
   $^{2}$Industrial Institute, Donetsk National Technical University, Pokrovsk, Ukraine\\
   $^{3}$Institute of Applied Mathematics and Mechanics, \\ National Academy of Sciences of Ukraine}

\begin{document}
\maketitle

\begin{abstract}                
A mechanical system consisting of a rigid body and attached Kirchhoff plates under the action of three independent controls torques is considered.
The equations of motion of such model are derived in the form of a system of coupled nonlinear ordinary and partial differential equations.
The operator form of this system is represented as an abstract differential equation in a Hilbert space.
A feedback control law is constructed such that the corresponding infinitesimal generator is dissipative.
\end{abstract}



\section{Introduction}
Problems of the aerospace industry and robotics stimulate the development of new methods for mathematical modeling and control design for complex mechanical systems with elastic elements.
In particular, it is a well-known fact that the vibrations of flexible parts of satellites influence significantly their dynamics, so that a rigid body model is not acceptable in stability and control investigations for such distributed parameter systems~(see, e.g.,~\cite{BM1968,LM1980,STY1986}).
This brings the motivation for studying the controllability and stabilization issues for infinite-dimensional mathematical models of flexible structures with strings, beams, and plates.
Without pretending to be complete, we refer to the monographs in this area~\cite{Krabs1992,Lag-Leug,DZ2006,M2013,Zuyev}.

The stabilization problem for a thin plate with boundary control has already received attention in~\cite{L1989,LT1991,HL1995,GZ2016}
A mathematical model of a rigid body with the Kirchhoff plate has been considered in the paper~\cite{Z2010}. It is assumed there that the body rotates around the fixed axis and its angular acceleration is taken as the control. The reachable sets for the linearized representation of such a system with modal coordinates have been analysed in~\cite{ZN2015}.

The purpose of our present paper is to derive a nonlinear model of a rotating rigid body with two Kirchhoff plates and propose a stabilizing feedback control for this model.
We will consider spatial rotations of the system and treat the three independent torques, applied to the body, as control inputs.
This framework is considered as a mathematical model of a satellite with solar panels controlled by jet thrusters or flywheels.


\section{Nonlinear model of the rotational motion}

Consider a mechanical system that consists of a rigid
body and two elastic plates~(Fig.~1).

\begin{figure}[ht!]
  \begin{center}
  \includegraphics[width=1.\linewidth]{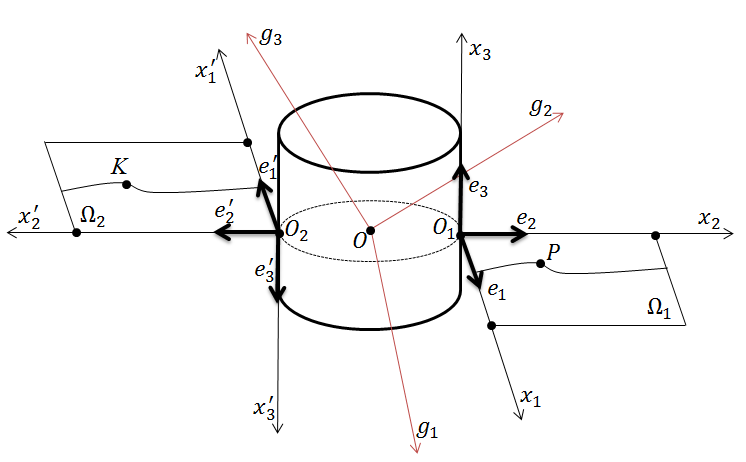}
   \caption{Rigid body with elastic plates.} \label{fig-1}
  \end{center}
\end{figure}

Let $(g_1,g_2,g_3)$ be the unit vectors of a fixed Cartesian frame.
Suppose that the rigid body is endowed with two Cartesian frames $O_1x_1x_2x_3$ and $O_2x'_1x'_2x'_3$,
and their basis vectors  $(e_1,e_2,e_3)$ and $(e'_1,e'_2,e'_3)$, respectively, are related as $e'_1=-e_1$,  $e'_2=-e_2$, and $e'_3=-e_3$.

Let $M=f_1e_1+f_2e_2+f_3e_3$ be the torque of external forces applied to the rigid body. \emph{We will treat the components $(f_1,f_2,f_3)^T\in\mathbb R^3$ as control inputs and consider the problem of defining a state feedback law in order to stabilize the moving frame $(e_1,e_2,e_3)$ in the direction of $(g_1,g_2,g_3)$ and to damp the vibrations of the plates}.
Note that a similar problem for an absolutely rigid body was solved in the book~\cite{Z1975}, and the problem of partial stabilization was considered in~\cite{Z2001,NDST2009,Z2006}.

In this paper, we assume that two rectangular plates are attached to the rigid body, so that in the undeformed state their median surfaces are located on the planes $O_1 x_1 x_2$ and $O_2 x_1' x_2'$, respectively.
At time $t$, the coordinates of a point $P$ on the median surface of the first plate can be represented in the frame $O_1 x_1 x_2 x_3$ as
$$
P=(x_1,x_2,w_1(x_1,x_2,t)),\;(x_1,x_2)\in \Omega_1=[0,l_1]\times[0,l_2].
$$
Similarly, the coordinates of a point $K$ on the median surface of the second plate are as follows (in the frame $O_2 x_1' x_2' x_3'$):
$$
K = (x'_1,x'_2,w_2(x'_1,x'_2,t)),\; (x'_1,x'_2)\in\Omega_2=[0,l'_1]\times[0,l'_2].
$$
Thus, the functions $w_1(x_1,x_2,t)$ and $w_2(x_1',x_2',t)$ define the transverse displacements for the case of small deformations of the plates.

In order to describe the motion of the considered mechanical system, we assume that the center of mass of the rigid body (point~$O$) is fixed and expand the vectors $O O_1$ and $O O_2$ with respect to the moving frames:
$$
O O_1 = d_1 e_1 + d_2 e_2 +d_3 e_3,\; O O_2 = d_1' e_1' + d_2' e_2' +d_3' e_3'.
$$
Then the absolute velocities of the points $P$ and $K$ are, respectively,
\begin{equation}\label{Z-N-4}
v_P=\omega
\times r_P+\dot{w_1}e_3
\end{equation}
and
\begin{equation}\label{Z-N-4-1}
v_K=\omega
\times r_K-\dot{w}_2e_3,
\end{equation}
where $\omega=\omega_1 e_1 + \omega_2 e_2 + \omega_3 e_3$ is the angular velocity vector of the rigid body,
\begin{equation}\label{Z-N-5}
r_P=(x_1+d_1)e_1+(x_2+d_2)e_2+(w_1+d_3)e_3,
\end{equation}
\begin{equation}\label{Z-N-5-1}
r_K=-(x'_1+d'_1)e_1-(x'_2+d'_2)e_2-(w_2+d'_3)e_3.
\end{equation}

Thus, from formulas~(\ref{Z-N-4}) and~(\ref{Z-N-4-1}) with taking into account~(\ref{Z-N-5}) and~(\ref{Z-N-5-1}), we get:
\begin{multline} \label{Z-N-6}
v_P=(\omega_2(d_3+w_1)-\omega_3(d_2+x_2))e_1+\\
+
(\omega_3(d_1+x_1)-\omega_1(d_3+w_1))e_2+\\
+(\omega_1(d_2+x_2)-\omega_2(d_1+x_1)+\dot{w}_1)e_3,
\end{multline}
\begin{multline} \label{Z-N-6-1}
v_K=-(\omega_2(d'_3+w_2)-\omega_3(d'_2+x'_2))e_1-\\
-
(\omega_3(d'_1+x'_1)-\omega_1(d'_3+w_2))e_2-\\
-(\omega_1(d'_2+x'_2)-\omega_2(d'_1+x'_1)+\dot{w}_2)e_3.
\end{multline}

According to the Kirchhoff plate model (cf.~\cite{L1989}), we write the following partial differential equations with respect to $w_1$ and $w_2$:
\begin{equation}\label{Z-N-1}
\begin{aligned}
\ddot{w}_1+a_1^2\left( \frac{\partial^2}{\partial x_1^2}+\frac{\partial^2}{\partial x_2^2} \right)^2 w_1&=(x_1+d_1)\dot{\omega}_2-(x_2+d_2)\dot{\omega}_1,\\
& \text{for}\;(x_1,x_2) \in \Omega_1,
\end{aligned}
\end{equation}
\begin{equation}\label{Z-N-1-1}
\begin{aligned}
\ddot{w}_2+a_2^2\left( \frac{\partial^2}{\partial x_1'^2}+\frac{\partial^2}{\partial x_2'^2} \right)^2w_2& =(x'_1+d'_1)\dot{\omega}_2-(x'_2+d'_2)\dot{\omega}_1,\\
& \text{for}\;(x_1',x_2') \in \Omega_2,
\end{aligned}
\end{equation}
where
$a_j^2=\frac{E_jh_j^3}{12\rho_j(1-\nu_j^2)}$ is the stiffness parameter of the $j$-th plate,
$E_j$ is Young's modulus, $\nu_j$ is Poisson's ratio,
$\rho_j$  is the area density, and
 $h_j$ is the thickness of the $j$-th plate.
The right-hand sides of~(\ref{Z-N-1}) and~\eqref{Z-N-1-1} contains the inertia forces because of the rotational motion of the rigid body (cf.~\cite{Lure}). In the differential equations~(\ref{Z-N-1}) and~\eqref{Z-N-1-1},
only linear terms with respect to the displacements, angular velocities, and derivatives of these quantities are taken into account (this is the linearized model of the plates' vibrations).

We assume that the plates are supported at their boundaries, which yields the following boundary conditions:
\begin{equation}\label{Z-N-2}
\begin{aligned}
& w_j|_{\partial\Omega_j}=0,\\
& \left.\frac{\partial^2w_j}{\partial n^2}\right|_{\partial\Omega_j}=0,\quad j=1,2.
\end{aligned}
\end{equation}
Here $\left.\frac{\partial w_j}{\partial n}\right|_{\partial\Omega_j}$ is the normal derivative of $w_j$ evaluated at the boundary of $\Omega_j$.

To derive the equations of motion of the rigid body-carrier, we exploit the angular momentum equation with respect to the fixed point $O$ (see, e.g.,~\cite{Lure}):
\begin{equation}\label{Z-N-3}
\dot{K}+\omega\times K=M,
\end{equation}
where $K=K_1 e_1 + K_2 e_2 + K_3 e_3$ is the angular momentum of the system, and $\dot K$ stands for the local derivative of $K$ in the moving frame $(e_1,e_2,e_3)$, i.e.
\begin{equation}
\dot K = \dot K_1 e_1 + \dot K_2 e_2 + \dot K_3 e_3.
\label{local}
\end{equation}

In the sequel, we use formulas~\eqref{Z-N-5}--\eqref{Z-N-6-1} to express the angular momentum $K$  for the mechanical system under consideration:
\begin{equation}
K = I\omega+  K_{p1}+ K_{p2},
\label{kinmoment}
\end{equation}
where $I$ is the tensor of inertia of the rigid body, and
$$
K_{p1} = \int\limits_{\Omega_1} r_P\times v_P\rho_1 dx,\; K_{p2}=\int\limits_{\Omega_2} r_K\times v_K \rho_2 dx'.
$$
We will assume that $(e_1,e_2,e_3)$ are the principal axes of inertia of the rigid body to simplify computations,
so that $I=diag(I_1,I_2,I_3)$.

Let us now compute the terms $K_{p1}$ and $K_{p2}$ for~\eqref{kinmoment}.
Formulas~(\ref{Z-N-5}) and (\ref{Z-N-6}) imply that,
\begin{equation}\label{Z-N-7}
K_{p1}=\int\limits_{\Omega_1} \rho_1(K_{11}e_1+K_{12}e_2+K_{13}e_3)dx,
\end{equation}
where
\begin{multline*}
K_{11}=\omega_1[(x_2+d_2)^2+(w_1+d_3)^2]-
\omega_2(x_1+d_1)(x_2+d_2)-\\-\omega_3(x_1+d_1)(w_1+d_3)+\dot{w}_1(x_2+d_2),\\
K_{12}=-\omega_1(x_1+d_1)(x_2+d_2)+\omega_2[(x_1+d_1)^2+(w_1+d_3)^2]-\qquad \quad\\
-\omega_3(x_2+d_2)(w_1+d_3)-\dot{w}_1(x_1+d_1),\\
K_{13}=-\omega_1(x_1+d_1)(w_1+d_3)-\omega_2(x_2+d_2)(w_1+d_3)+\qquad \qquad \qquad\\
+\omega_3[(x_1+d_1)^2+(x_2+d_2)^2].
\end{multline*}
Similarly, from formulas (\ref{Z-N-5-1}) and (\ref{Z-N-6-1}) it follows that
\begin{equation}\label{Z-N-71}
K_{p2}=\int\limits_{\Omega_2} \rho_2(K_{21}e_1+K_{22}e_2+K_{23}e_3)dx',
\end{equation}
where
\begin{multline*}
K_{21}=\omega_1[(x'_2+d'_2)^2+(w_2+d'_3)^2]-
\omega_2(x'_1+d'_1)(x'_2+d'_2)-\\-\omega_3(x'_1+d'_1)(w_2+d'_3)+\dot{w}_2(x'_2+d'_2),\\
K_{22}=-\omega_1(x'_1+d'_1)(x'_2+d'_2)+\omega_2[(x'_1+d'_1)^2+(w_2+d'_3)^2]-\qquad \quad\\
-\omega_3(x'_2+d'_2)(w_2+d'_3)-\dot{w}_2(x'_1+d'_1),\\
K_{23}=-\omega_1(x'_1+d'_1)(w_2+d'_3)-\omega_2(x'_2+d'_2)(w_2+d'_3)+\qquad \qquad \qquad\\
+\omega_3[(x'_1+d'_1)^2+(x'_2+d'_2)^2].
\end{multline*}

Then the angular momentum~\eqref{kinmoment} with the use of~\eqref{Z-N-7} and~\eqref{Z-N-71} can be rewritten as follows:
\begin{equation}
\begin{aligned}
K=&(I+J)\omega+\sum_{n=1}^{2}\int\limits_{\Omega_n}\dot{w}_n\left(x_{2}+d_{2n}\right) \rho_n dx \,e_1-\\
&-
\sum_{n=1}^{2}\rho_n\int\limits_{\Omega_n}\dot{w}_n\left(x_{1}+d_{1n}\right)dx \, e_2+ K_\delta,
\end{aligned}\label{Z-N-8}
\end{equation}
where $J=(J_{ij})$ is the tensor of inertia for the mechanical system with ``frozen'' plates (i.e. when the plates are considered as rigid bodies),
\begin{align*}
J_{11}=\sum_{n=1}^{2}\rho_n\int\limits_{\Omega_n}\left((x_{2}+d_{2n})^2+d_{3n}^2\right)dx, \\
J_{12}=J_{21}=-\sum_{n=1}^{2}\rho_n\int\limits_{\Omega_n}(x_{1}+d_{1n})(x_{2}+d_{2n})dx, \\
J_{22}=\sum_{n=1}^{2}\rho_n\int\limits_{\Omega_n}\left((x_{1}+d_{1n})^2+d_{3n}^2\right)dx, \\
J_{23}=J_{32}=-\sum_{n=1}^{2}\rho_n d_{3n}\int\limits_{\Omega_n}(x_{2}+d_{2n})dx, \\
J_{33}=\sum_{n=1}^{2}\rho_n\int\limits_{\Omega_n}\left((x_{1}+d_{1n})^2+(x_{2}+d_{2n})^2\right)dx,\\
J_{31}=J_{13}=-\sum_{n=1}^{2}\rho_n d_{3n}\int\limits_{\Omega_n} (x_{1}+d_{1n})dx,
\end{align*}
and the term $K_\delta$ is of order $O\Bigl(\|\omega\|\left(\|\dot w_1\|+\|\dot w_2\|\right)\Bigr)$ for small $\dot w_j$.
By computing the local derivative (in the sense of~\eqref{local}) for the angular momentum $K$ given by formula~\eqref{Z-N-8}, we get:
$$
\dot K_1  = (J_{11}+I_1)\dot{\omega}_1+J_{12}\dot{\omega}_2
+J_{13}\dot{\omega}_3+\sum_{n=1}^{2}\rho_n\int\limits_{\Omega_n}\ddot{w}_n(x_{2}+d_{2n})dx,
$$
$$
\dot K_2  =J_{21}\dot{\omega}_1+(J_{22}+I_2)\dot{\omega}_2
+J_{23}\dot{\omega}_3-\sum_{n=1}^{2}\rho_n\int\limits_{\Omega_n}\ddot{w}_n(x_{1}+d_{1n})dx,
$$
\begin{equation}
\dot K_3  = J_{31}\dot{\omega}_1+J_{32}\dot{\omega}_2+(J_{33}+I_3)\dot{\omega}_3,
\label{dotK}
\end{equation}
where the nonlinear terms with respect to the derivatives of $\omega_j$ and $w_n$ are omitted.

We expand the cross product $\omega\times K= \sum_{j=1}^3 (\omega\times K)_j e_j$ to see that
{\small
\begin{multline}\nonumber
(\omega\times K)_1 =
 \omega_2\left(J_{31}\omega_1+J_{32}\omega_2+(J_{33}+I_3)\omega_3\right)-J_{21}\omega_1 \omega_3-\\-
\omega_3\left((J_{22}+I_2)\omega_2+J_{23}\omega_3-\sum_{n=1}^{2}\rho_n\int\limits_{\Omega_n}\dot{w}_n\left(x_{1}+d_{1n}\right) dx\right),
\end{multline}
\begin{multline}\nonumber
(\omega\times K)_2 =-
\omega_1\left(J_{31}\omega_1+J_{32}\omega_2+(J_{33}+I_3)\omega_3\right)+J_{12}\omega_2\omega_3+\\
+\omega_3\left(\sum_{n=1}^{2}\rho_n\int\limits_{\Omega_n}\dot{w}_n\left(x_{2}+d_{2n}\right)dx+
(J_{11}+I_1)\omega_1+J_{13}\omega_3\right),
\end{multline}
\begin{multline}\label{Z-N-10}
(\omega\times K)_3 =
\omega_1\left(J_{21}\omega_1+(J_{22}+I_2)\omega_2+J_{23}\omega_3\right)-\\
-\omega_2\left((J_{11}+I_1)\omega_1+J_{12}\omega_2+J_{13}\omega_3
\right)-\\
- \sum_{n=1}^{2}\rho_n\int\limits_{\Omega_n}
\dot{w}_n\bigl(\omega_2\left(x_{2}+d_{2n}\right)+\omega_1 \left(x_{1}+d_{1n}\right)\bigr)dx.
\end{multline}
}

By putting together  the formulas (\ref{dotK}) and (\ref{Z-N-10}) and expressing the values of $\ddot{w}_1$ and $\ddot{w}_2$ from~(\ref{Z-N-1}) and~(\ref{Z-N-1-1}), we write the components of the angular momentum equation~(\ref{Z-N-3}) as follows:
{\begin{multline}\label{Z-N-11}
\left(I_1+\sum_{n=1}^{2}\rho_nd_{3n}^2l_{1n}l_{2n}\right)\dot{\omega}_1+J_{13}\dot{\omega}_3=\\
=f_1+
\omega_3\left[J_{21}\omega_1+(J_{22}+I_2)\omega_2+J_{23}\omega_3\right]-\\
-\omega_2\left[J_{31}\omega_1+J_{32}\omega_2+(J_{33}+I_3)\omega_3\right]+
\\
+\sum_{n=1}^{2}\rho_n\int\limits_{\Omega_n}\{a_n^2(x_{2}+d_{2n})\Delta^2w_n-\omega_3\dot{w}_n(x_{1}+d_{1n})\}dx,\\
 \left(I_2+\sum_{n=1}^{2}\rho_n d_{3n}^2l_{1n}l_{2n}\right)\dot{\omega}_2+J_{23}\dot{\omega}_3=\\
 =f_2+
\omega_1\left[J_{31}\omega_1+J_{32}\omega_2+(J_{33}+I_3)\omega_3\right]-\\
-\omega_3\left[(J_{11}+I_1)\omega_1+J_{12}\omega_2+J_{13}\omega_3\right]-\\
-
\sum_{n=1}^{2}\rho_n\int\limits_{\Omega_n}\{a_n^2(x_{1}+d_{1n})\Delta^2w_n+\omega_3\dot{w}_n(x_{2}+d_{2n})\}dx,\\
J_{31}\dot{\omega}_1+J_{32}\dot{\omega}_2+(J_{33}+I_3)\dot{\omega}_3=\\
=f_3-
\omega_1\left[J_{21}\omega_1+(J_{22}+I_2)\omega_2+J_{23}\omega_3\right]+\\
+\omega_2\left[(J_{11}+I_1)\omega_1+J_{12}\omega_2+J_{13}\omega_3\right]+\\
+\sum_{n=1}^{2}\rho_n\int\limits_{\Omega_n}\{(x_{1}+d_{1n})\omega_1+(x_{2}+d_{2n})\omega_2\}\dot{w}_ndx,
\end{multline}}
where the nonlinear terms with respect to the derivatives are omitted.

In order to  rewrite the above differential equations in the normal form with respect to $\dot \omega_j$, we compute the inverse matrix $J^{-1}= \hat J =  ({\hat J}_{ij})$:
$$
{\hat J}_{11} = \frac{(I_2 + \sum_{n=1}^{2}\rho_n d_{3n}^2 l_{1n} l_{2n})(I_3+J_{33})-J_{23}^2}{D},
$$
$$
{\hat J}_{12} = {\hat J}_{21} =\frac{J_{13} J_{23}}{D},
$$
$$
{\hat J}_{13} ={\hat J}_{31} =-\frac{(I_2 + \sum_{n=1}^{2}\rho_n d_{3n}^2 l_{1n} l_{2n})J_{13}}{D},
$$
\begin{equation}
{\hat J}_{22} = \frac{(I_1 + \sum_{n=1}^{2}\rho_n d_{3n}^2 l_{1n} l_{2n})(I_3+J_{33})-J_{13}^2}{D},
\label{Jinv}
\end{equation}
$$
{\hat J}_{23} ={\hat J}_{32} =-\frac{(I_1 + \sum_{n=1}^{2}\rho_n d_{3n}^2 l_{1n} l_{2n})J_{23}}{D},
$$
$$
{\hat J}_{33} = \frac{(I_1 + \sum_{n=1}^{2}\rho_n d_{3n}^2 l_{1n} l_{2n})(I_2 + \sum_{n=1}^{2}\rho_n d_{3n}^2 l_{1n} l_{2n})}{D},
$$
{\small
\begin{multline*}
D = (I_1 + \sum_{n=1}^{2}\rho_n d_{3n}^2 l_{1n} l_{2n})(I_2 + \sum_{n=1}^{2}\rho_n d_{3n}^2 l_{1n} l_{2n})(I_3+J_{33})-\\-J_{13}^2 (I_2 + \sum_{n=1}^{2}\rho_n d_{3n}^2 l_{1n} l_{2n})-J_{23}^2 (I_1 + \sum_{n=1}^{2}\rho_n d_{3n}^2 l_{1n} l_{2n}).
\end{multline*}
}
Note that the denominator $D$ in formulas~\eqref{Jinv} is strictly positive at least for sufficiently small moments of inertia of the plates $J_{ik}$ compared with the moments of inertia of the carrier body $I_i$.
In particular, this condition is satisfied for sufficiently small area densities $\rho_j$ (i.e., for sufficiently thin plates). Thus, we assume that $D\neq 0$ in the sequel. Then the differential equations~\eqref{Z-N-11} can be written in the form
\begin{equation}
\begin{pmatrix}\dot\omega_1 \\ \dot\omega_2 \\ \dot\omega_3 \end{pmatrix} = \hat J \begin{pmatrix}\phi_1 \\ \phi_2 \\ \phi_3 \end{pmatrix},
\label{omega_DE}
\end{equation}
where $\phi_i$ denotes the right-hand side of the $i$-th equation in~\eqref{Z-N-11}.

We write the Poisson kinematic equations to ensure the condition that the frame $(g_1,g_2,g_3)$ is fixed in the inertial space:
\begin{equation}\label{Z-N-12}
\dot{g}_i=-\omega\times g_i, \qquad i=\overline{1,3}.
\end{equation}
Let $g_i=g_{i1}e_1+g_{i2}e_2+g_{i3}e_3$, then system~(\ref{Z-N-12}) takes the following coordinate form:
\begin{equation}\label{Z-N-13}
\begin{split}
&\dot{g}_{i1}=\omega_3g_{i2}-\omega_2g_{i3}, \\
&\dot{g}_{i2}=\omega_1g_{i3}-\omega_3g_{i1}, \\
&\dot{g}_{i3}=\omega_2g_{i1}-\omega_1g_{i2}, \qquad i=\overline{1,3}.
\end{split}
\end{equation}

For the Cartesian frames $(g_1,g_2,g_3)$ and $(e_1,e_2,e_3)$ of the same orientation, the system of
differential equations
(\ref{Z-N-1}), (\ref{Z-N-1-1}), (\ref{Z-N-2}), (\ref{Z-N-11}), (\ref{Z-N-13})
has the following particular solution with $f_1=f_2=f_3=0$:
\begin{equation}\label{Z-N-14}
w(x,t)=0, \quad \omega_i(t)=0, \quad g_{ij}(t)=\delta_{ij},  \qquad i,j=\overline{1,3},
\end{equation}
where $\delta_{ij}$ is the Kronecker symbol.

To study the stabilization problem for the equilibrium~(\ref{Z-N-14}),
we introduce new variables $\widetilde{g}_{ij}(t)=g_{ij}(t)-\delta_{ij}$ and consider the equations
of perturbed motion for~(\ref{Z-N-13}):
\begin{equation}\label{Z-N-15}
\begin{aligned}
&\dot{\widetilde{g}}_{11}=\omega_3\widetilde{g}_{12}-\omega_2\widetilde{g}_{13}, \quad
\dot{\widetilde{g}}_{12}=\omega_1\widetilde{g}_{13}-\omega_3(\widetilde{g}_{11}+1), \\
&\dot{\widetilde{g}}_{13}=\omega_2(\widetilde{g}_{11}+1)-\omega_1\widetilde{g}_{12},\\
&\dot{\widetilde{g}}_{21}=\omega_3(\widetilde{g}_{22}+1)-\omega_2\widetilde{g}_{23}, \quad
\dot{\widetilde{g}}_{22}=\omega_1\widetilde{g}_{23}-\omega_3\widetilde{g}_{21},\\
&\dot{\widetilde{g}}_{23}=\omega_2\widetilde{g}_{21}-\omega_1(\widetilde{g}_{22}+1),\\
&\dot{\widetilde{g}}_{31}=\omega_3\widetilde{g}_{32}-\omega_2(\widetilde{g}_{33}+1),\;
\dot{\widetilde{g}}_{32}=\omega_1(\widetilde{g}_{33}+1)-\omega_3\widetilde{g}_{31},\\
&\dot{\widetilde{g}}_{33}=\omega_2\widetilde{g}_{31}-\omega_1\widetilde{g}_{32}.
\end{aligned}
\end{equation}

We consider a modified energy functional
\begin{equation}\label{Z-N-16}
\cal E =T+U+\displaystyle\frac12\sum\limits_{i,j=1}^3\alpha_i\widetilde{g}_{ij}^2
\end{equation}
with positive parameters $\alpha_i$, where
$$T=\frac12\left (I_1\omega_1^2+I_2\omega_2^2+I_3\omega_3^2+\int\limits_{\Omega_1} v_P^2 \rho_1dx+~\int\limits_{\Omega_2} v_K^2 \rho_2 dx \right)$$ is the
kinetic energy of the system, and
$$U~=~\displaystyle\frac12 \sum_{n=1}^{2}\int\limits_{\Omega_n}(\triangle w_n(x,t))^2 a_n^2 \rho_n  dx$$ is the potential energy of elastic deformations according to the Kirchhoff model.
Here $\Delta = \frac{\partial^2}{\partial x_1^2}+ \frac{\partial^2}{\partial x_2^2}$ is the Laplace operator.
For future use, we introduce the Lyapunov functional $V$ as a quadratic approximation of $\cal E$:
{\small
$$
2V=(I_1+J_{11})\omega_1^2+(I_2+J_{22})\omega_2^2+(I_3+J_{33})\omega_3^2+2J_{12}\omega_1\omega_2 +
$$
$$
+
2J_{13}\omega_1\omega_3
+2 J_{23}\omega_2\omega_3+\sum_{i,j=1}^3 \alpha_i {\tilde g_{ij}}^2+\sum_{n=1}^{2}\rho_n \int\limits_{\Omega_n}{\dot w}_n^2 dx+
$$
\begin{equation}
+
\sum_{n=1}^{2}\rho_n \int\limits_{\Omega_n}\left\{2\dot w_n [\omega_1(d_{2n}+x_{2})-\omega_2(d_{1n}+x_{1})]+a_n^2 (\Delta w_n)^2\right\}dx.
\label{Lyapunov}
\end{equation}
}

Let us compute the time derivative of the functional~(\ref{Z-N-16}) along the trajectories of~(\ref{Z-N-1}), (\ref{Z-N-11}), (\ref{Z-N-15}):
$$
\dot V = \left(\dot K_1 - \alpha_2 \tilde g_{23}+\alpha_3 \tilde g_{32}\right)\omega_1 + \left(\dot K_2 + \alpha_1 \tilde g_{13}-\alpha_3 \tilde g_{31}\right)\omega_2 +
$$
$$
+\left(\dot K_3 + \alpha_2 \tilde g_{21}-\alpha_1 \tilde g_{12}\right)\omega_3+
$$
\begin{equation}
+\sum_{n=1}^{2} \int\limits_{\Omega_n}  \left\{\Delta w_n \Delta\dot w_n - \dot w_n\Delta^2 w_n\right\}a_n^2 \rho_n dx,
\label{Vdot1}
\end{equation}
where the expressions for $\dot K_i$ are given by~\eqref{dotK}.
If the partial derivatives of $w_i(x,t)$ of the fourth order in $x$ and the first order in $t$ are continuous and the boundary conditions~\eqref{Z-N-2} are satisfied, then the integration by parts in formula~\eqref{Vdot1} leads to the following identities:
$$
\int_{\Omega_n} \left\{\Delta w_n \Delta\dot w_n - \dot w_n\Delta^2 w_n\right\}dx =0,\quad n=1,2.
$$
Using these identities and expressing $\dot K_i$ from the equation~\eqref{Z-N-3}, we rewrite formula~\eqref{Vdot1} as
$$
\dot V = \left(f_1 - (\omega\times K)_1 - \alpha_2 \tilde g_{23}+\alpha_3 \tilde g_{32}\right)\omega_1 +
$$
$$
+\left(f_2 - (\omega\times K)_2 + \alpha_1 \tilde g_{13}-\alpha_3 \tilde g_{31}\right)\omega_2 +
$$
$$
+\left(f_3 - (\omega\times K)_3 + \alpha_2 \tilde g_{21}-\alpha_1 \tilde g_{12}\right)\omega_3,
$$
where $(\omega\times K)_i$ are given by~\eqref{Z-N-10}.

{\em
To stabilize the trivial solution of the system~\eqref{Z-N-1},\eqref{Z-N-1-1}, \eqref{Z-N-2}, \eqref{Z-N-11}, \eqref{Z-N-15}, we define a feedback control from the condition
\begin{equation}
\dot{V} = -k(\omega_1^2+\omega_2^2+\omega_3^2)\le 0,
\label{Vdot_condition}
\end{equation}
where $k$ is a positive constant}. It is easy to see that condition~\eqref{Vdot_condition} corresponds to the following choice of controls:
$$
f_1 = -k\omega_1 +(\omega\times K)_1 + \alpha_2 \tilde g_{23}-\alpha_3 \tilde g_{32},
$$
$$
f_2 = -k\omega_2 +(\omega\times K)_2 -\alpha_1 \tilde g_{13}+\alpha_3 \tilde g_{31},
$$
\begin{equation}
f_3 = -k\omega_3 +(\omega\times K)_3+\alpha_1 \tilde g_{12}-\alpha_2 \tilde g_{21}.
\label{feedback}
\end{equation}
Note that the time-derivative~\eqref{Vdot_condition} is not negative definite, and the finite-dimensional method described in~\cite{Z2016} is not directly applicable
to establish strictly decreasing behavior of $V$ along the trajectories of the closed-loop system.

\section{Operator form of the dynamical equations}

To formulate our main result, we rewrite the equations of motion of the mechanical system under consideration in the operator form. We introduce the real linear space $H={\mathring H}^2(\Omega_1)\times{\mathring H}^2(\Omega_2)\times L^2(\Omega_1)\times L^2(\Omega_2)\times \mathbb R^{12}$, whose elements are denoted as
\begin{multline*}
\xi=\begin{pmatrix} u_1 \\ v_1 \\ u_2 \\ v_2 \\ \omega \\ \tilde g \end{pmatrix}:
u_n\in {\mathring H}^2(\Omega_n),\; v_n\in L^2(\Omega_n),\; n=1,2,\\ \omega = \begin{pmatrix} \omega_1 \\ \omega_2 \\ \omega_3\end{pmatrix}\in\mathbb R^3,\; \tilde g = \begin{pmatrix}\tilde g_{11} \\ \tilde g_{12} \\ \vdots \\ \tilde g_{33}\end{pmatrix}\in\mathbb R^9.
\end{multline*}
Here ${\mathring H}^2(\Omega_n)$ is the Sobolev space of the functions $u\in H^2(\Omega_n)$ having zero trace on $\partial \Omega_n$.
The inner product of
$$
\xi^1=\begin{pmatrix} u_1^1 \\ v_1^1 \\ u_2^1 \\ v_2^1 \\ \omega^1 \\ {\tilde g}^1 \end{pmatrix}\in H\;\;\text{and}\;\;
\xi^2=\begin{pmatrix} u_1^2 \\ v_1^2 \\ u_2^2 \\ v_2^2 \\ \omega^2 \\ {\tilde g}^2 \end{pmatrix}\in H
$$
is defined as
$$
\left<\xi^1,\xi^2\right>_H = \sum_{n=1}^{2}\rho_n \int\limits_{\Omega_n}\Bigl\{a_n^2 \Delta u_n^1(x)  \Delta u_n^2(x) + v_n^1(x) v_n^2(x) +
$$
$$
+ (\omega^2_1 v_n^1(x)+\omega^1_1 v_n^2(x))(d_{2n}+x_{2})-$$
$$
-(\omega^2_2 v_n^1(x)+ \omega^1_2 v_n^2(x))(d_{1n}+x_{1})
\Bigr\}dx+
$$
\begin{equation}
+\Bigl((I+J)\omega^1,\omega^2\Bigr)+\sum_{i,j=1}^3\alpha_i {\tilde g}^1_{ij} {\tilde g}^2_{ij}.
\label{scalar_p}
\end{equation}
Using the Cauchy--Schwarz and Friedrichs' inequalities, it can be shown that the norm $\|\xi\|_H=\sqrt{\left<\xi,\xi\right>_H}$ is equivalent to the standard norm in ${\mathring H}^2(\Omega_n)\times L^2(\Omega_n)\times \mathbb R^{12}$.
Thus, $\left(H, \left<\cdot,\cdot\right>_H\right)$ is a Hilbert space.

We define an unbounded operator $A: D(A)\to H$ and a bounded linear operator $B:\mathbb R^3\to H$ in the following way:
\begin{equation}
A: \xi=\begin{pmatrix} u_1 \\ v_1 \\u_2 \\ v_2 \\ \omega \\ \tilde{g} \end{pmatrix} \mapsto
A\xi = \begin{pmatrix} u_1^\xi \\ v_1^\xi \\  u_2^\xi \\ v_2^\xi \\ \omega^\xi \\ {\tilde{g}}^\xi \end{pmatrix}\in H,
\label{op_A}
\end{equation}
\begin{equation}
B: f=\begin{pmatrix} f_1 \\ f_2 \\ f_3 \end{pmatrix}\mapsto Bf = \begin{pmatrix} u_1^f \\ v_1^f \\ u_2^f \\ v_2^f \\ \omega^f \\ {\tilde{g}}^f \end{pmatrix}\in H,
\label{op_B}
\end{equation}
where
$$
\omega^\xi_i = ({\hat J}_{i1}\omega_3-\hat J_{i3}\omega_1)\left[J_{21}\omega_1+(J_{22}+I_2)\omega_2+J_{23}\omega_3\right]+
$$
$$
+({\hat J}_{i2}\omega_1-\hat J_{i1}\omega_2)\left[J_{31}\omega_1+J_{32}\omega_2+(J_{33}+I_3)\omega_3\right]+
$$
$$
+({\hat J}_{i3}\omega_2-\hat J_{i2}\omega_3)\left[(J_{11}+I_1)\omega_1+J_{12}\omega_2+J_{13}\omega_3\right]+
$$
$$
+\sum_{n=1}^{2}\rho_n \int\limits_{\Omega_n}\left(\hat J_{i3}[(x_{1}+d_{1n})\omega_1 + (x_{2}+d_{2n})\omega_2]\right. -
$$
$$
\left.
-[\hat J_{i1}(x_{1}+d_{1n})+\hat J_{i2}(x_{2}+d_{2n})]\omega_3\right) v_n(x)dx+
$$
$$
+\sum_{n=1}^{2}\rho_n a_n^2\int\limits_{\Omega_n}  \left(\hat J_{i1} (x_{2}+d_{2n})-\hat J_{i2} (x_{1}+d_{1n})\right)\Delta^2 u_n(x)dx,
$$
$$
v_n^\xi(x) = -a_n^2\Delta^2 u_n(x)+(x_{1}+d_{1n})\omega^\xi_2 - (x_{2}+d_{2n})\omega^\xi_1,
$$
$$
u_n^\xi(x) = v_n(x),\quad n=1,2,
$$
$$
{\tilde{g}}^\xi_{11} = \omega_3\widetilde{g}_{12}-\omega_2\widetilde{g}_{13},\;
{\tilde{g}}^\xi_{12} =\omega_1\widetilde{g}_{13}-\omega_3(\widetilde{g}_{11}+1),
$$
$$
{\tilde{g}}^\xi_{13} =\omega_2(\widetilde{g}_{11}+1)-\omega_1\widetilde{g}_{12},
$$
$$
{\tilde{g}}^\xi_{21} =\omega_3(\widetilde{g}_{22}+1)-\omega_2\widetilde{g}_{23},\;
{\tilde{g}}^\xi_{22} =\omega_1\widetilde{g}_{23}-\omega_3\widetilde{g}_{21},
$$
$$
{\tilde{g}}^\xi_{23} =\omega_2\widetilde{g}_{21}-\omega_1(\widetilde{g}_{22}+1),
$$
$$
{\tilde{g}}^\xi_{31} =\omega_3\widetilde{g}_{32}-\omega_2(\widetilde{g}_{33}+1),\;
{\tilde{g}}^\xi_{32} =\omega_1(\widetilde{g}_{33}+1)-\omega_3\widetilde{g}_{31},
$$
$$
{\tilde{g}}^\xi_{33} =\omega_2\widetilde{g}_{31}-\omega_1\widetilde{g}_{32},
$$
$$
v_n^f(x) =\sum_{k=1}^3\left(\hat J_{2k}(x_{1}+d_{1n})-\hat J_{1k}(x_{2}+d_{2n})\right)f_k,
$$
$$
u_n^f(x)=0,\;{\tilde{g}}^f =0,\;
\omega^f_i = \sum_{k=1}^3\hat J_{ik}f_k,\quad i=1,2,3,
$$
and the coefficients $\hat J_{ik}$ are given in~\eqref{Jinv}.
The domain of definition of the nonlinear operator $A$ has the form
\begin{equation}
\label{Z-N-22}
D(A)=\Bigg\{\xi\in H  \,\left| \; \begin{array}{l}u_n\in {\mathring H}^4(\Omega_n), \; v_n \in {\mathring H}^2(\Omega_n),
\\ \left.\frac{\partial^2u_n}{\partial x_1^2}\right|_{\partial \Omega_n}=\left.\frac{\partial^2u_n}{\partial x_2^2}\right|_{\partial _n}=0, n=1,2\end{array}\right.\Bigg\}.
\end{equation}

\section{Main result}

We consider a nonlinear control system governed by the following abstract differential equation in $H$:
\begin{equation}
\frac{d}{dt}\xi(t) = A\xi(t) + B f,
\label{abstract_DE}
\end{equation}
where $\xi(t)\in H$ is the state state, $f\in \mathbb R^3$ is the control, and the operators $A$ and $B$ are given by~\eqref{op_A} and~\eqref{op_B}.
If the functions $w_1(x,t)$, $w_2(x,t)$, $\omega(t)$, $\tilde g(t)$ define a classical solution of  system~\eqref{Z-N-1}, \eqref{Z-N-1-1}, \eqref{Z-N-2}, \eqref{omega_DE}, \eqref{Z-N-15} with a control  $f=f(t)$ on the interval $t\in {\cal I}=[t_0,T)$, $T\le +\infty$,
then by direct substitution we verify that the corresponding function
\begin{equation}
\xi(t) = \begin{pmatrix}w_1(\cdot,t) \\ \dot w_1 (\cdot,t) \\w_2(\cdot,t) \\ \dot w_2 (\cdot,t) \\ \omega(t) \\ {\tilde g}(t)\end{pmatrix}\in D(A)\subset H
\label{xi_form}
\end{equation}
satisfies the equation~\eqref{abstract_DE} with $f=f(t)$ on $t\in \cal I$.
Thus, {\em we consider the differential equation~\eqref{abstract_DE} as a mathematical model of the considered mechanical system}.

Let us represent the feedback control~\eqref{feedback} by an operator $G:H\to \mathbb R^3$ defined on the state space of system~\eqref{abstract_DE}:
\begin{equation}
G: \xi=\begin{pmatrix}u_1 \\ v_1 \\ u_2 \\ v_2 \\ \omega \\ \tilde g \end{pmatrix}\mapsto f=Gx=\begin{pmatrix}f_1^\xi \\ f_2^\xi \\ f_3^\xi \end{pmatrix},
\label{op_G}
\end{equation}
{
$$
f_1^\xi = -
\omega_3\left(J_{21}\omega_1+(J_{22}+I_2)\omega_2+J_{23}\omega_3\right)-
$$
$$
-k\omega_1 +\alpha_2 \tilde g_{23}-\alpha_3 \tilde g_{32} + \omega_2\left(J_{31}\omega_1+J_{32}\omega_2+(J_{33}+I_3)\omega_3\right)+
$$
$$
+
\omega_3\sum_{n=1}^{2}\rho_n \int\limits_{\Omega_n} \left(x_{1}+d_{1n}\right)v_n(x) dx,
$$
$$
f_2^\xi =
-
\omega_1\left(J_{31}\omega_1+J_{32}\omega_2+(J_{33}+I_3)\omega_3\right)
-
$$
$$
-k\omega_2  -\alpha_1 \tilde g_{13}+\alpha_3 \tilde g_{31}
+\omega_3\left(
(J_{11}+I_1)\omega_1+J_{12}\omega_2+J_{13}\omega_3\right)+
$$
$$
+\omega_3\sum_{n=1}^{2}\rho_n \int\limits_{\Omega_n} \left(x_{2}+d_{2n}\right)v_n(x)dx,
$$
$$
f_3^\xi = -\omega_2\left((J_{11}+I_1)\omega_1+J_{12}\omega_2+J_{13}\omega_3
\right)-
$$
$$
-k\omega_3 +\alpha_1 \tilde g_{12}-\alpha_2 \tilde g_{21}+\omega_1\left(J_{21}\omega_1+(J_{22}+I_2)\omega_2+J_{23}\omega_3\right)-
$$
$$
-\sum_{n=1}^{2}\rho_n \int\limits_{\Omega_n}
[\left(x_{1}+d_{1n}\right)\omega_1+\left(x_{2}+d_{2n}\right)\omega_2]v_n(x)dx,
$$
}
where $k$ and $\alpha_i$ are arbitrary positive constants.
Then the closed-loop system~\eqref{abstract_DE} with  $f=G\xi$ takes the form
\begin{equation}
\frac{d}{dt}\xi(t)= F\xi(t),\quad F=A+BG,
\label{closed_loop}
\end{equation}
where the domain of definition of the unbounded nonlinear operator $F:D(F)\to H$ is dense in $H$, $D(F)=D(A)$.

As it was noted above, the classical solutions of system~\eqref{Z-N-1}, \eqref{Z-N-1-1}, \eqref{Z-N-2}, \eqref{omega_DE}, \eqref{Z-N-15} correspond to the functions $\xi(t)\in D(A)$ according to the rule~\eqref{xi_form}.
Then the condition $\dot V\le 0$ for the time derivative of the functional $V$ along the trajectories of the closed-loop system can be rewritten as
\begin{equation}
\left<F\xi,\xi\right>_H\le 0
\label{dis_ineq}
\end{equation}
for the corresponding element $\xi\in D(F)=D(A)$. This is a consequence of the definition of the functional $V$ in~\eqref{Lyapunov} and the inner product $\left<\cdot,\cdot\right>_H$ in~\eqref{scalar_p}.
Thus, inequality~\eqref{dis_ineq} implies the following result concerning the operator $F$ of the closed-loop system~\eqref{closed_loop}.

{\bf Theorem~1.}\label{thm_dis}
{\em The operator $F:D(F)\to H$ is dissipative and $\overline{D(F)}=H$.}

We denote by $I_H$ the identity operator on $H$.
If $F$ is closed and the image of  $I_H-\lambda F$  coincides with $H$ for $\lambda>0$, then Theorem~1 implies that $F$  is the infinitesimal generator of a strongly continuous semigroup of nonlinear operators
$\{S(t)\}_{t\geq0}$ on $H$
because of the Crandall--Liggett theorem~(cf.~\cite{Barbu}).
 Then the mild solution of the Cauchy problem for~\eqref{closed_loop}
with the initial condition $\xi(0)=\xi^0$ is defined by the formula
\begin{equation}
\xi(t)=S(t)\xi^0,\quad t\ge 0,
\label{mild}
\end{equation}
for any $\xi^0\in H$. {\em Under these assumptions, the trivial solution of the abstract differential equation~\eqref{closed_loop} is stable in the sense of Lyapunov because of the dissipativity inequality $\left<F\xi,\xi\right>_H\le 0$ ($\dot V\le 0$)}.
Note that the mild solutions given by formula~\eqref{mild} are classical if $\xi^0\in D(F)$.

\section{Numerical simulations}
In order to illustrate the transient behaviour of the proposed controller, we perform a numerical simulation of the closed-loop system~\eqref{closed_loop}.
For this purpose we consider the case of identical rectangular domains $\Omega_1=\Omega_2=[0,l_1]\times [0,l_2]$ and introduce finite-dimensional approximations of the displacements
$w_j(x_1,x_2,t)=q_j(t)W_1(x_1)W_2(x_2)$, $j=1,2$, for equations~\eqref{Z-N-1} and~\eqref{Z-N-1-1}, respectively.
Here $W_j(x)=\sin\left(\frac{\pi m_j x}{l_j}\right)$, $x\in [0,l_j]$, $m_j\in \mathbb N$, are taken as eigenfunctions of the Sturm--Liouville problem with the boundary conditions
$W_j(0)=W_j(l_j)=0$, $j=1,2$. Let us consider the first flexible mode only ($m_1=m_2=1$) and apply the Ritz--Galerkin method (cf.~\cite{ZS2015}) for the nonlinear closed-loop system~\eqref{Z-N-1}--\eqref{Z-N-2},~\eqref{omega_DE},~\eqref{Z-N-15},~\eqref{feedback} to derive its finite-dimensional approximation of the form
\begin{equation}
\dot X (t) = \Phi(X(t)),\quad X(t) \in {\mathbb R}^{16},
\label{finite_dim}
\end{equation}
whose state vector is
$$
X = \left(\tilde g_{11}, \tilde g_{12}, ..., \tilde g_{33} , \omega_1,\omega_2,\omega_3, q_1,q_2,\dot q_1,\dot q_2 \right).
$$
We choose the following initial data and parameters for the simulation (the dimensions of physical quantities are omitted to simplify notations):
\begin{equation}
X(0)=\left(0,0,0, 0,-1,1, 0,-1,-1, 0, ..., 0\right),
\label{IC}
\end{equation}
$$
\begin{aligned}
&
l_1=\frac{l_2}{2}=1,\; d_1=d_1'=0,\; d_2=d_2'=1,\; d_3=d_3'=0, \\
& \rho_1=\rho_2=1,\; a_1=a_2 = \frac{1}{2},\; J=I,\; \alpha_1 = \alpha_2=\alpha_3=1.
\end{aligned}
$$
The above $X(0)$ corresponds to an equilibrium of the considered mechanical system with $g_1=e_1$, $g_2 = e_3$, and $g_3=-e_2$.
In this case, the considered stabilization problem (steering the closed-loop system to its trivial equilibrium) means rotation about the $x_1$-axis by the angle $\pi/2$
with simultaneous damping of the vibration modes.

  \setlength{\unitlength}{1cm}
\begin{figure}[ht!]
  \begin{center}
  \includegraphics[width=1.\linewidth]{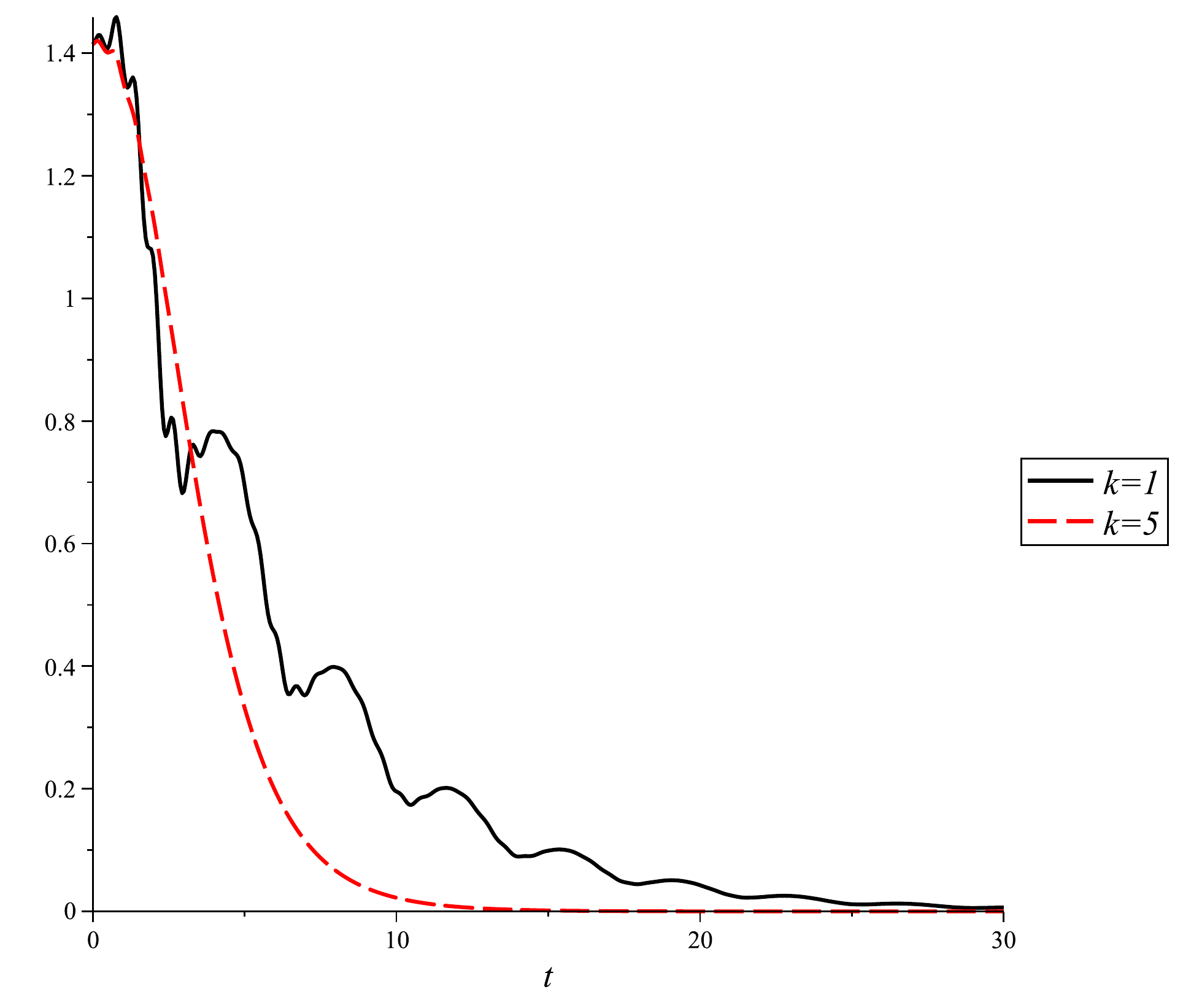}
  \put(-10,9){$\|X(t)\|$}
   \caption{Euclidean norm $\|X(t)\|$ of the solution of~\eqref{finite_dim},~\eqref{IC}.} \label{fig-2}
  \end{center}
\end{figure}

Fig.~\ref{fig-2} shows the behavior of solutions to the Cauchy problem~\eqref{finite_dim},~\eqref{IC} for different values of the tuning parameter $k>0$ appearing in the feedback law~\eqref{feedback}.
We see that the proposed controller can be used to steer the state of system~\eqref{finite_dim} to zero as $t\to +\infty$, and the higher value of $k$ ($k=5$) leads to faster convergence of $X(t)$ to $X=0$.

\section{Conclusions}

In this paper, a new mathematical model of a controlled mechanical system consisting of a rotating body and elastic Kirchhoff plates has been introduced.
The model is described by the system of nonlinear ordinary and partial differential equations~\eqref{Z-N-1}--\eqref{Z-N-2},~\eqref{omega_DE},~\eqref{Z-N-15}, or, equivalently,
by the abstract differential equation~\eqref{abstract_DE} in the Hilbert space $H$.
A state feedback control has been proposed explicitly in the form~\eqref{feedback} to ensure that the time derivative of a Lyapunov functional is non-positive.

The main theoretical contribution of this work establishes the dissipativity property for the infinitesimal generator $F$ in~\eqref{closed_loop}.
Although the results of numerical simulations illustrate the efficiency of the proposed controller,
the question about {\em asymptotic stability} (or even {\em partial asymptotic stability} in the sense of~\cite{Z2003}) remains open for the infinite-dimensional closed-loop dynamics.


\end{document}